\documentclass[a4paper]{article}
\usepackage[T1]{fontenc}
\usepackage{multicol,verbatim,graphicx,epsfig,amssymb,amsmath,mathrsfs}
\usepackage{amsfonts}
\usepackage{amssymb}
\usepackage{latexsym}
\usepackage{graphicx}
\usepackage{epstopdf} 

\newtheorem{theorem}{Theorem}[section]
\newtheorem{corollary}[theorem]{Corollary}
\newtheorem{lemma}[theorem]{Lemma}
\newtheorem{definition}[theorem]{Definition}
\newtheorem{remark}[theorem]{Remark}
\newtheorem{example}[theorem]{Example}

\newcommand{\proof}{\medskip \noindent {\textbf{Proof.} \ \ }}
\newcommand{\qed}{\null\hfill $\Box\;\;$ \medskip}

\begin{document}
\parbox{1mm}

\begin{center}
{\bf {\sc \Large On fields inspired with the polar $HSV-RGB$ theory of Colour}
\footnote{\it Mathematics Subject
Classification (2000): \sl 92C55, 
%Biomedical imaging and signal processing
 12E30. 
% 	Field arithmetic,
\newline {\it Key words and phrases:}
RGB representation of colour,
polar colour space,
parabolic complex numbers
HSV theory
\newline {\it Acknowledgement.} This paper was
supported by Grant VEGA 2/0178/14 and by the Slovak-Ukrainian joint research project "Vector valued measures and integration in polarized vector spaces".}
}
\end{center}

\vskip 12pt

\begin{center}
{\bf J\'an Halu\v{s}ka\footnote{
J\'an Halu\v{s}ka, Mathematical Institute, Slovak Academy of Sciences, Gre\v{s}\'{a}kova~6, 040~00 Ko\v{s}ice, Slovakia, e-mail:~jhaluska@saske.sk}
}\end{center}

\begin{comment}
\subsection*{Abstract}
\footnotesize

A three-polar, cf.~\cite{Gregor-Haluska1}, $HSV-RGB$ Colour space $\triangle$ was introduced and studied. It was equipped with operations of addition, subtraction, multiplication, and (partially) division. Achromatic Grey Hues  form an ideal $\mathfrak{S}$.  Factorizing  $\triangle$  by the ideal $\mathfrak{S}$,  we obtain a field $\triangle | \mathfrak{S}$. An element (i.e an individual  Colour) in  $\triangle | \mathfrak{S}$ is a triplet of three triangular coefficients. The set of all triangular coefficients is a subset of a semi-field of  parabolic-complex functions. For the parabolic-complex number set, cf.~\cite{Harkin}.
\end{comment}

\normalsize
%%%%%%%%%%%%%%%%%%%%%%%%%
\section{Introduction}
\subsection{Traditional Colour  theories}
%%%%%%%%%%%%%%%%%%%%%%%%%
There are more different sources of the theory of Colour which approach to the subject from different sides and are  complementary in this sense.
 
\paragraph{Computer graphics}
The HSV (Hue-Saturation-Value) theory is the most common representation of points in an RGB (Red-Green-Blue) color \textit{technical model}. Computer graphics pioneers developed the HSV model in the 1970s for computer graphics applications (A. R. Smith in 1978, also in the same issue, A. Joblove and H. Greenberg). A HSV theory is used today in color pickers, in image editing software, and less commonly in image analysis and computer vision. A rather extensive explanation of the present State of Art in industry we can find in \cite{So}.

\paragraph{Biophysics}
 Th. Young and H.~Helmholtz proposed a trichromatic theory. Their theory states that the human retina contains dispersed photo-sensitive clusters, where each of these clusters consists of three types of sensitive cones which peaks in \textit{short} (420--440~nm), \textit{middle} (530--540~nm), and long (560--580~nm) wavelengths. Weighting a total light power spectrum by the individual spectral sensitivities of the three types of cone cells 
gives three effective stimulus values; these three values make up a tristimulus specification of the objective color of the light spectrum. In Fig.~\ref{birds}, there are schematic behaviours of these sensitive clusters of cells.

\paragraph{Fine arts}
 For the revelatory theories of Colour written by authors who are not mathematicians (artistic photographers, visual artists), c.f.~\cite{Briggs,Hirsch,Hunt}.

%%%%%%%%%%%%%%%%%%%%%%%%%%%%%%%%%%%%%%%%%%%%
\subsection{Terminology}
%%%%%%%%%%%%%%%%%%%%%%%%%%%%%%%%%%%%%%%%%%%%
For terminology, basic and also advanced concepts about Colour, we refer to~\cite{Stockman}, Chapter~11 (Vol. III, Vision and vision optics; Chap. 11, Color vision mechanism).

%%%%%%%%%%%%%%%%%%%%%%% %%%%%%%%%%%%%%%%%%%%%%%%%%%%%%
\subsection{Comments to modelling of Colour}
%%%%%%%%%%%%%%%%%%%%%%%%%%%%%%%%%%%%%%%%%%%%%%%%%%%%%%
 
A  reflected  electro-magnetic vibration energy  is  filtered  with the  (human)  tristimulus apparatus  in the  eye retina into three  functions (also called the $RGB$ stimulus curves). Here is an information loss, because we see in the wavelength interval approximately  350 nm -- 750 nm although there is some power output theoretically within the whole $[0, + \infty)$. Vibrations within other vibration intervals are partially perceived by other senses (hearing, touch). Also variuos animals have various intervals within  they can see.   

The vision process continues in the brain  where the  obtained three stimuli curves are aggregated back. The result is a registry of Colour of the object. 
In our theory, this aggregation is a linear combination of three basic colours (poles) with the coefficients which are complex functions defined on the interval $[0, +\infty)$ (= all possible frequencies).   
 
\begin{figure}
\centering
   \includegraphics[totalheight=8cm]{%obrazky/
   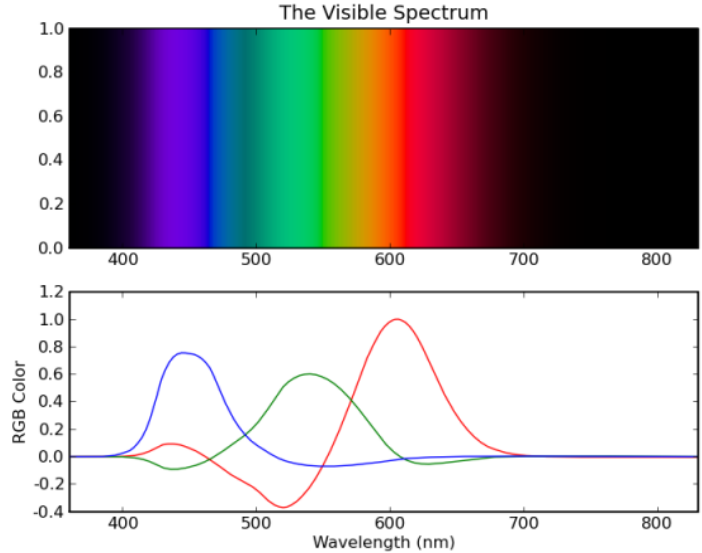}
  \caption{A realistically shaped tristimulus reaction on a Colour} %\cite[p.~4]{XXXXX}.}
  \label{birds} % Figure 1
\end{figure}

In Fig.~\ref{birds}, we see  that  realistic stimulus curves may have parts which are particularly in the negative (under the $x$-axis), i.e., the accession to the resulting Colour may happen also such that the parts of curves "absorbs" energy. 
Some practical aspects about $RGB$ in computer graphics we can find in \cite{Pascale}. 

%%%%%%%%%%%%%%%%%%%%%%%%%%%%%%%%%%%%%%%%%

%%%%%%%%%%%%%%%%%%%%%%%%%%%%%%%%%%%%%%%%%

%%%%%%%%%%%%%%%%%%%%%%%%%%%%%%%%%%%%%%%%%%%%%%%%%%%%%%%%%
\begin{figure}[h]
 \centering
 \includegraphics[totalheight=5cm]{%obrazky/
 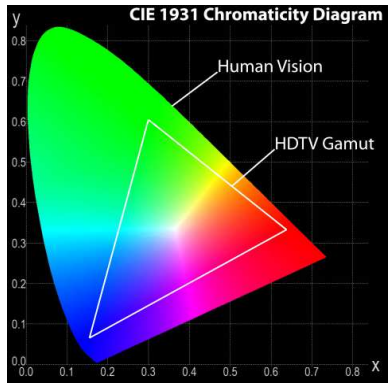}
 \caption{A horse-shoe planar cut of the visible colour space; Black = non visible}
 \label{cie1931}
\end{figure}

%%%%%%%%%%%%%%%%%%%%%%%%%%%%%%%%%%%%%%%%%%%%%%%%%%%%%%%%

%%%%%%%%%%%%%%%%%%%%%%%%%%%

\begin{figure}[h]
\centering
 \includegraphics[totalheight=4cm]{%obrazky/
 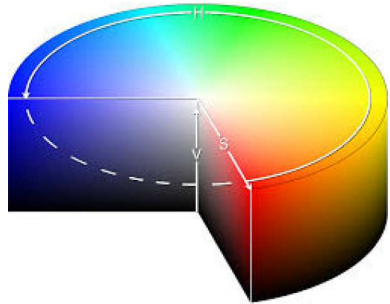}
 \caption{A Colour wheel with colours of the  equal energy (the upper surface)} \label{wheelsatu}
\end{figure}

%%%%%%%%%%%%%%%%%%%%%%%%%%%%%%%%%%%%%%%%%%%%%%%%%%%%%
Outside a  certain area of the domain, the colour aggregated from  $R, B, G$ curves  is  not very exact  because of distorted perceptions by human senses, see Fig.~\ref{cie1931}. 
 
We use an imagination of the Colour Hues in the Colour wheel. The reason is that there are colour hues which are not in the linear rainbow palette (e.g., Pink).
We suppose that angles  of the basic $RGB$ colours in the colour wheel correspond to the angles  $0 \equiv 2 \pi$, $2 \pi/3$, $4\pi/3$, respectively, which is  approximately true in reality.

  An incoming white sun light is Colours is decomposed into a sum of three $RGB$ curves. So it is also in our model. But abstracting from the natural perception,  the number three for basic colours  is  not substantial. Our polar theory  works for arbitrary natural $K\geq 3$.  For a review, in the realm of animals, there are mono-chromatic   Arctic mammals; most of mammals have sensibility only for two colours, they are dichromats; birds and insects are mostly tetrachromats. Concerning primates, the human vision is trichromatic. The record keeps the Mantis shrimp's  vision with $K=12$, cf. \cite{mantis}. 

There are also  artificial colour schemes coming and used in the industry for some good reasons. E.g., we know the  CMYK (cyan, magenta yellow, black), RGYB (red, green, yellow, blue) systems, and others. 

%%%%%%%%%%%%%%%%%%%%%%%%%%%%%%%%%%%%%%%%%
\subsection{The point and interval characteristics of light}
%The set of all trianglular coefficients as a proper subset of parabolic-complex functions is a semi-field}
%%%%%%%%%%%%%%%%%%%%%%%%%%%%%%%%%%%%%%%%%
Hue is a point characteristic (it is determined in a point), the  saturation and brightness are  "interval" characteristics, i.e., they are determined for an interval, not at a point (similarly as  the notions concavity-convexity has no sense at a points).

\begin{figure}[h]
 \centering
 \includegraphics[totalheight=4cm]{%obrazky/
 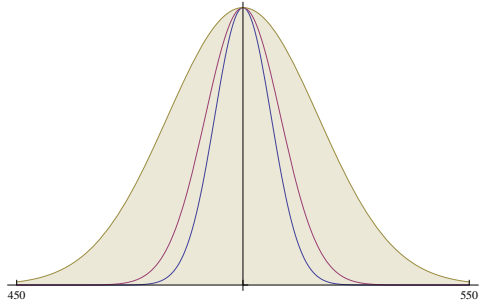}
 \caption{Monotonicity of Saturation; Brightness and Hue are fixed}
 \label{obr3}
\end{figure}

\begin{figure}[h]
 \centering
 \includegraphics[totalheight=4cm]{%obrazky/
 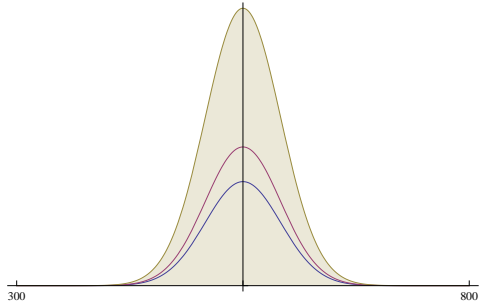}
 \caption{Monotonicity of Brightness; Saturation and Hue are fixed}
 \label{obr2}
\end{figure}

\begin{figure}[h]
 \centering
 \includegraphics[totalheight=4cm]{%obrazky/
 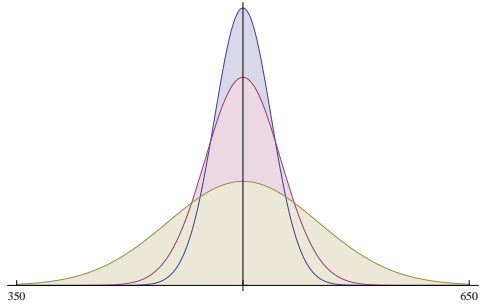}
 \caption{Non-comparable Saturation and Brightness, Hue is fixed}
 \label{obr1}
\end{figure}

\paragraph{Hue} of Colour is the wavelength within the visible-light spectrum at which the energy output from a source is greatest. This is shown as the peak of the sum of the three intensity curves in the accompanying graphs of intensity versus wavelength. In the illustrative examples in the pictures Fig.\ref{obr3},  Fig.\ref{obr2}, Fig.\ref{obr1}, all $3\times 3$ = nine  colors there have the  same  Hue with a wavelength 500~nm, in the yellow-green portion of the spectrum.

\paragraph{Saturation} is an expression for the relative bandwidth of the visible output from a light source. In the diagram~Fig.\ref{obr3}, the notion of saturation is represented relatively, by the steepness of the slopes of the curves. Here, the blue  curve represents a color having the greatest saturation. As saturation increases, the color with the same Hue appears more "pure." As saturation decreases, colors appear more "washed-out."

\paragraph{Brightness} is a relative expression of the intensity of the energy output of a visible light source. It can be expressed as a total energy value (different for each of the curves in the diagram Fig.~\ref{obr2}), or as the amplitude at the wavelength where the intensity is greatest. Energy is imagined as the area under the curve.  In the picture the blue curve has the lowest brightness.

As we can see, Saturation and Brightness  are generally non-comparable parameters of Colour, cf. Fig.\ref{obr1}.

One commonly supposes that  all possible colours can be specified  according to these three parameters and that Colors can be represented in terms of the $RGB$ components. Thus the whole information in the  $RBG$ Colour theory is contained  the three tristimulus curves.  A concept of triangular  coefficients   mathematically reflects this idea.

%%%%%%%%%%%%%%%%%%%%%%%%%%%%%%%%%%%%%%%%%%%%%
\subsection{Semi-field of triangular coefficients}\label{coef}
%%%%%%%%%%%%%%%%%%%%%%%%%%%%%%%%%%%%%%%%%%%%%

A \textit{semi-field} $\mathbb{X}$ is a set equipped with an algebraic structure with binary operations of addition $(+)$
and multiplication $(\cdot)$, where $(\mathbb{X},+)$ is a commutative semi-group. $(X, \cdot, 1)$ is a multiplicative group with the unit $1$, and multiplication is distributive with respect to addition from both sides. For a review of semi-fields, c.f.~\cite{Vechto}.
A semi-field $\mathbb{X}$ is called a \textit{semi-field with zero} if there exists an element $\mathbb{O}\in \mathbb{X}$ such that 
$\mathbf{x} \cdot \mathbb{O} = \mathbb{O}$ and $\mathbf{x} + \mathbb{O} = \mathbf{x}$ for every $\mathbf{x} \in \mathbb{X}$ and the distributivity of multiplication from both sides is preserved for the extended system.

\begin{example}\sl
Let $\mathbb{R}_{+,0} := (0,\infty) \cup \{0\} = [0, \infty)$ be a ray with all structures heredited from the real line $\mathbb{R}$. This is one of \textit{trivial real semi-field with zero}.    
\end{example}

\begin{definition} \rm  We say that a set 
$$ T  = \{ q + \psi(f) \varepsilon \mid q \geq 0, \psi(f) \in \mathbb{R}^{[0, +\infty)}, f \in [0, + \infty) \}$$
is called to be  the \textit{set of all  triangular coefficients}, where $\varepsilon$ is the  \textit{parabolic imaginary unit}, $|\varepsilon|=1, \varepsilon^2 = 0$, cf.~\cite{Harkin}. 

If   $[a + b(f) \varepsilon] \in T$ and $[A + B(f) \varepsilon] \in T$  are  two  triangular coefficients,   we define the operations of addition, multiplication, and division as  following.
For every $ a, A \geq 0$; $b(f), B(f) \in \mathbb{R}^{[0, +\infty)}$, $f \in [0, + \infty)$,
\begin{eqnarray}\label{tplus}
 [a+  b(f)\varepsilon] +  [A +  B(f)\varepsilon] :=  [a + A ] + [b(f)+B(f)]\varepsilon, \end{eqnarray}
 \begin{eqnarray}\label{trnumbermulti}  [a + b(f)\varepsilon] \cdot  [A + B(f)\varepsilon] 
:=  a A + [a B(f) + A b(f)]\varepsilon, \end{eqnarray}
\begin{eqnarray}\label{trnumber2} A \neq 0 \implies  \frac{a + b(f)\varepsilon}{A + B(f)\varepsilon} :=   \frac{a}{A} + \frac{A b(f) - aB(f)}{A^2}\varepsilon.  
 \end{eqnarray}
\end{definition}

For the better reading, triangular coefficients are written in square brackets in the sequel of the paper.

\begin{lemma} The set $T$ with the operations defined in (\ref{tplus}), (\ref{trnumbermulti}), (\ref{trnumber2}) is a semi-field.
\end{lemma}
  \proof Saving the denotation of the previous definition we see that
$[a+  b(f)\varepsilon] +  [A +  B(f)\varepsilon] \in T$, $[a + b(f)\varepsilon] \cdot  [A + B(f)\varepsilon] \in T$, and  for $A \neq 0$, $ \frac{a + b(f)\varepsilon}{A + B(f)\varepsilon} \in T$.

The assertion is an obvious  enlargement of the result for parabolic-complex  numbers to functions 
defined on the non-negative real ray. \qed

In our polar theory of Colour, a role of scalars will play elements of the semi-field $T$.
 
\begin{corollary}
For every function $B(f), f  \in [0, + \infty)$, and by (\ref{trnumber2}),
 $$  A \neq 0 \implies \frac{1}{A + B(f)\varepsilon} =    \frac{1}{A} + 
 \frac{ - B(f)}{A^2} \varepsilon \in T. $$ 
\end{corollary}
%%%%%%%%%%%%%%%%%%%%%%%%%%%%%%%%%%%%%%%%%%%%%%%%%%%%%%%%%
%%%%%%%%%%%%%%%%%%%%%%%%%%%%%%%%%%%%%%%%%%%%%%%%%%%%%%%%

%%%%%%%%%%%%%%%%%%%%%%%%%%%%%%%%%%%%%%%%%%%%%%%%%%%%%%%%%
%%%%%%%%%%%%%%%%%%%%%%%%%%%%%%%%%%%%%%%%%%%%%%%%%%%%%%%%

%%%%%%%%%%%%%%%%%%%%%%%%%%%%%%%%%%%%%%%%%%%%%%%%%%%%%%%%%
%%%%%%%%%%%%%%%%%%%%%%%%%%%%%%%%%%%%%%%%%%%%%%%%%%%%%%%%

\begin{comment}
\begin{remark} \rm
The parabolic-complex semi-field is substantially used in our Col\-our theory  when we deal with multiplication/division of colours. For addition/subtraction of Colours all three (elliptic, hyperbolic, or parabolic) types of complex numbers are equally suitable.

The reason of using parabolic-complex numbers (for multiplication of Colours in the sequel) is that the set of all parabolic numbers (note that the whole set of all parabolic complex numbers is not a field) has  a proper two-dimensional subset  which is a non-trivial semi-field.   The standard (i.e., elliptic, Gaussian) complex numbers are not suitable since it has proper subsets which  are only one-dimensional semi-fields. The hyperbolic generalized complex numbers are not suitable for dealing with  achromatic hues in the sequel. For the satisfactory complete explanation of geometry an review of the generalized complex numbers,  cf.~\cite{Harkin}.  
\end{remark}
\end{comment}
%%%%%%%%%%%%%%%%%%%%%%%%%%%%%%%%%%%%%%%%%%%%%%%%%%%%%%%%% 

\subsection*{The aim of the paper}\textit{ 
To model a  Colour space  as a three-polar  field. }

%%%%%%%%%%%%%%%%%%%%%%%%%%%%%%%%%%%%%%%%%%%%%%%%%%%%%%%%%%%%%%%%%%%%%%%%%%%%%%%%%%%%%%%%%%%%%%%
\section{Mastering polar Colour spaces}
%%%%%%%%%%%%%%%%%%%%%%%%%%%%%%%%%%%%%%%%%%%%%%%%%%%%%%%%%%%%%%%%%%%%%%%%%%%%%%%%%%%%%%%%%%%%%%% 
\subsection{Poles and the  definition of  Colour space} \label{pooles} 
%%%%%%%%%%%%%%%%%%%%%%%%%%%%%%%%%%%%%%%%%%%%%%%%%%%%%%%%%%%%%%%%%%%%%%%%%%%%%%%%%%%%%%%%%%%%%%% 

The real line is two-polar and we are using the obvious polar operators, i.e., the signs $(+)$ and $(-)$. In this sense, we will also understand three signs (poles)  $R, G, B$ in this paper. The poles $R, G, B$ could be understand also as the \textit{generalized signs}.

In fact, poles  can be chosen as objects  of various mathematical nature. In this paper, we use the following 
finite set of poles:
%%%%%%%%%%%%%%%%%%%
 $$\label{Examplesx} \mathbb{A} = \left\{ R= 1 + 0 \imath, G = -\frac{1}{2} + \frac{\sqrt{3}}{2} \imath, B= -\frac{1}{2}  -\frac{\sqrt{3}}{2}\imath \right\}, $$  a set of three vertices of an equilateral triangle in the  elliptic complex plane, where $\imath$ is a usual (elliptic) complex unit, $\imath^2 = -1$, $|\imath| = 1$. 
Sometimes we will equivalently speak that poles are \textit{polar operators}. 

\begin{remark}\rm 
Which properties are asked from poles in general? Various other objects of different nature  can also be used as poles. E.g., functions, operators, the set $\mathbb{A}$ when replacing complex unit $\imath$ by the complex units $\varepsilon, \kappa$ of the parabolic ($\varepsilon^2=0$) or hyperbolic ($\kappa^2=1$) complex units, respectively, etc. 
In our case, we use the  operators $R, G, B \in \mathbb{C}$ which  
\begin{enumerate}
\item  
  are applicable to "all objects" (similarly as signs plus and minus);
  \item  fulfil  the  condition  $R + G + B = 0$ (presence of the \textit{white point} of the colour space);
  \item $R,G,B$ are different and non-collinear points;
  \item    a symmetry  in some sense of the set $\mathbb{A}$ is  desirable.
  \end{enumerate}
  All these terms will be precised below. 
\end{remark}

\begin{definition} \rm
Let us denote by
\begin{equation} 
\begin{array} {rcl}\triangle & := &\{ R[r+\rho(f)\varepsilon]+ G[g + \gamma(f)\varepsilon] + B[b + \beta(f)\varepsilon] \\ && \mid [r +\rho(f)\varepsilon] \in T,  [g + \gamma(f)\varepsilon] \in T, [b + \beta(f)\varepsilon] \in T\},\end{array} \label{3polarity}\end{equation}
where  $r, g, b  \geq 0$; $\rho(f), \gamma(f), \beta(f) 
\in \mathbb{R}^{[0, +\infty)}, f \in [0, +\infty)$.
The set $\triangle$ is called the \textit{colour space}.  An element of the colour space is called to be the \textit{colour}.
 \end{definition}

\begin{remark} \rm
(1) Although the functions $\rho(f), \gamma(f), \beta(f) \in \mathbb{R}^{[0, \infty)}$  are arbitrary in our paper, in the  praxis they are supposed to have "good" properties, e.g., they are supposed to be unimodal, continuous, etc. 
(2) Elements of $\triangle$ are (mixed elliptic--parabolic) bicomplex numbers. For (elliptic--elliptic) bicomplex numbers, cf. \cite{Shapiro}.
\end{remark}

%%%%%%%%%%%%%%%%%%%%%%%%%%%%%%%%%%%%%%%%%%%%%
\subsection {Achromatic part of Colour}\label{achroma}
%%%%%%%%%%%%%%%%%%%%%%%%%%%%%%%%%%%%%%%%%%%%%
Every Colour is a composite of chromatic (pure colours) and achromatic (grey and noises) parts. Both parts of colour can be derived  from the decomposition of the tristimulus sum into three  individual curves. In this paper, the achromatic part of colour is derived from  an element  $[a + \alpha(f)] \in T$, where $a \in [0,\infty)$ determines a Hue of Grey.\footnote {The so called Value, which is a term overtaken from the HSV Colour theory; not very apt for a mathematical theory.}  The function $\alpha(f) \in \mathbb{R}^{[0,\infty)}$ represents noise. We  incorporate an achromatic part of Colour into our theory 
 dealing with the following subsets (cuts) of 
the colour space $\triangle$:
\begin{definition} \rm  Let $[a + \alpha(f)\varepsilon] \in T$ and 
 $[r + \rho(f)\varepsilon] \in T$,  $[g + \gamma(f)\varepsilon] \in T$, $[b + \beta(f)\varepsilon] \in T$.
Denote by  
$$\mathfrak{o}_{[a + \alpha(f)\varepsilon]} :=  R[a +\alpha(f)\varepsilon]+ G[a + \alpha(f)\varepsilon] + B[a + \alpha(f))\varepsilon],  $$
$$\mathfrak{O} := \bigcup_{[a +  \alpha(f)\varepsilon] \in T} \mathfrak{o}_{[a + \alpha(f)\varepsilon]}   
\subset \triangle,$$
and
$$\mathfrak{s}_a := 
\{\mathbf{x} \in \triangle \mid \mathbf{x} = R[a +\rho(f)\varepsilon]+ G[a + \gamma(f)\varepsilon] + B[a + \beta(f)\varepsilon]\}, $$
$$\mathfrak{S} := \bigcup_{a \in [0, \infty)} \mathfrak{s}_{a},$$
where $a\geq 0$ and $\alpha(f), \rho(f), \gamma(f), \beta(f)$ are arbitrary  functions in $\mathbb{R}^{[0,\infty)}$. 
\end{definition} 

\begin{lemma} \label{congruence} Let $\mathbf{x}\ \in \triangle$. Let $\lambda \in \mathfrak{O}$. Then
$$ \mathbf{x} = \mathbf{x} + \lambda.$$
 \end{lemma}
\proof
From definition of $\mathfrak{O}$ it follows that the  triangular coefficients  $[r + \rho(f)\varepsilon]$, 
$[g + \gamma(f)\varepsilon]$, $[b + \beta(f)\varepsilon]$ in  $\mathbf{x} \in \triangle$ are ambiguous  since for every arbitrary $[a +  \alpha(f)\varepsilon] \in T$, there holds  
$$\begin{array}{rcl} \mathbf{x} & = \phantom{ + } & R [r + \rho(f)\varepsilon]  + G [g + \gamma(f)\varepsilon] + B [b + \beta(f)\varepsilon] \\
&  = \phantom{ + }  & R\{[r + \rho(f)\varepsilon] + [a + \alpha(f)\varepsilon]\} \\ 
&   \phantom { = } +  & G\{[(g + \gamma(f)\varepsilon] +[a + \alpha(f)\varepsilon]\} \\
& \phantom { = } + & B\{[(b + \beta(f)\varepsilon] +[a + \alpha(f)\varepsilon]\} 
 =\mathbf{x} + \lambda, \lambda \in \mathfrak{O}. \hskip1cm\Box
 \end{array}  $$
 
In the view of this lemma we introduce the following notion. 
\begin{definition}\rm
Let $X \in \mathbb{A}$. We say that colour $\mathbf{x}  \in \triangle$ is $X$-\textit{polarized} if 
it can be expressed in the form $\mathbf{x} = X [x +\xi(f) \varepsilon]$, $[x +\xi(f) \varepsilon] \in T$.
\end{definition}

\begin{definition}\rm
The  congruence given in Lemma~\ref{congruence} is  known also as the \textit{Cancellation law}, c.f. \cite{Gregor-Haluska1}).
So, (this holds concerning all arithmetic operations  in $\triangle$ we will define in the sequel of this paper). So,  we operate with the accuracy of congruent triples of  triangular coefficients.
\end{definition}

Physically Cancellation law means an ambiguity with respect to the achromatic parts of Colour. This is expressed with using of the phrase "with respect to Cancellation law". For the sake of simplicity and without loss of precision, this phrase will be often omitted in the text. 
 
%%%%%%%%%%%%%%%%%%%%%%%%%%%%%%%%%%%%%%%%%%%%%%%
\section{Arithmetic operations in $\triangle$}
%%%%%%%%%%%%%%%%%%%%%%%%%%%%%%%%%%%%%%%%%%%%%%%
Let us denote for the following sections: $$\begin{array}{rcl}\mathbf{x}& =& \{ R[r + \rho(f)\varepsilon] + G [g + \gamma(f)\varepsilon] + B [b + \beta(f)\varepsilon]\} \in \triangle, \\  \mathbf{y} & = & \{R [u + \sigma(f)\varepsilon] + G[v + \chi(f)\varepsilon] + B[t + \xi(f)\varepsilon]\} \in \triangle \end{array}$$ 
be two colours.

%%%%%%%%%%%%%%%%%%%%%%%%%%%%%%%%%%%%%%%%%%%%%%%%%%%%%%%%%%%%%%%%%
\subsection{Addition in $\triangle$ (\textit{Mixing of Colours})} 
%%%%%%%%%%%%%%%%%%%%%%%%%%%%%%%%%%%%%%%%%%%%%%%%%%%%%%%%%%%%%%%%%
We define:   
$$\begin{array}{rcl}\mathbf{x} \oplus \mathbf{y}& := & \phantom{+} R\{[r + \rho(f)\varepsilon] + [u + \sigma(f)\varepsilon]\} \\   & &+ G\{[g +  \gamma(f)\varepsilon] + [v + \chi(f)\varepsilon] \} \\ & & +  B\{[b+  \beta(f)\varepsilon] + [t + \xi(f)\varepsilon]\},\\ 
  & = & \phantom{+} R[\{r + u\} + \{ \rho(f) + \sigma(f) \}\varepsilon] \\   & &+ G[\{g + v\} + \{ \gamma(f) + \chi(f)\}\varepsilon]  \\ & & +  B[\{b+  t\} + \{\beta(f) + \xi(f)\}\varepsilon].
\end{array}$$

\begin{remark} Remind, that the result of operation of addition is with respect to Cancellation law, i.e.,
for every $\lambda_1, \lambda_2, \lambda_3 \in \mathfrak{O},$
 $$ \mathbf{x} \oplus \mathbf{y} = (\mathbf{x} + \lambda_1 )\oplus (\mathbf{y} + \lambda_2) = (\mathbf{x} \oplus \mathbf{y}) + \lambda_3.$$ 
\end{remark}

%%%%%%%%%%%%%%%%%%%%%%%%%%%%%%%%%%%%%%%%%%%%%%%%%%%%%%%%%%%%%%%%%%%%%%%
\subsection{Subtraction in $\triangle$ (\textit{Inverse colours})}
%%%%%%%%%%%%%%%%%%%%%%%%%%%%%%%%%%%%%%%%%%%%%%%%%%%%%%%%%%%%%%%%%%%%%%%

We define Subtraction in $\triangle$ as the addition of inverse elements, $$\mathbf{x} \ominus \mathbf{y} := \mathbf{x} \oplus (\ominus \mathbf{y}).$$ 
The inverse elements of the basic colours are defined as follows:
 $$\begin{array}{ccl}
 \ominus R [r +  \rho(f)\varepsilon] & := &G [r +  \rho(f)\varepsilon] + B[r +  \rho(f)\varepsilon],\\
 \ominus G [g +  \gamma(f)\varepsilon] & := &R [g +  \gamma(f)\varepsilon] + B[g +  \gamma(f)\varepsilon],\\
 \ominus B [b+  \beta(f)\varepsilon] & := &R [b+  \beta(f)\varepsilon] + G[b+  \beta(f)\varepsilon].
 \end{array} $$

 So, $$ \begin{array} {rcl}\ominus\mathbf{x} & := & \phantom{ + } R \{ [g +  \gamma(f)\varepsilon] + [b+  \beta(f)\varepsilon]\} \\ 
 && + G\{[r +  \rho(f)\varepsilon] + [b+  \beta(f)\varepsilon]\} \\ &&  + B\{ [r +  \rho(f)\varepsilon] + [g +  \gamma(f)\varepsilon] \}. \end{array}$$

And the subtraction is defined then as following:
$$ \begin{array}{ccl} \mathbf{x}\ominus \mathbf{y} &=& \phantom{\ominus}
\{R[r +  \rho(f)\varepsilon] + [g +  \gamma(f)\varepsilon] G+ B[b+  \beta(f)\varepsilon]\}  \\ &  & \ominus \{R[u + \sigma(f)\varepsilon] + G[v + \chi(f)\varepsilon] + B[t + \xi(f)\varepsilon]\} 
\\  &&\\  & := & \phantom{\oplus}
 \{R[r +  \rho(f)\varepsilon] + G[g +  \gamma(f)\varepsilon] + B[b+  \beta(f)\varepsilon]\} \\
&  & \oplus \{ G[u + \sigma(f)\varepsilon] + B[u + \sigma(f)\varepsilon]\} \\ & &  \oplus \{ R[v + \chi(f)\varepsilon] +B[v + \chi(f)\varepsilon] \} \\ &   & \oplus \{ R[t + \xi(f)\varepsilon] + G[t + \xi(f)\varepsilon]\}
\\ && \\
& =  & \phantom{ + }R\{ [r +  \rho(f)\varepsilon] + [v + \chi(f)\varepsilon] +[t + \xi(f)\varepsilon] \}  
\\ &  & + G\{[g +  \gamma(f)\varepsilon] +[u + \sigma(f)\varepsilon] + [t + \xi(f)\varepsilon] \} 
\\ & & + B\{[b+  \beta(f)\varepsilon] +[u + \sigma(f)\varepsilon] + [v + \chi(f)\varepsilon] \}.
\end{array} $$

\begin{remark}\rm
We can  replace polar operators $R,G,B$ with their \textit{inverse operators} $$C := \ominus R =-1,  M:= \ominus G = 1/2 - (\sqrt{3}/2))\imath, Y := \ominus B = 1/2 + (\sqrt{3}/2))\imath.$$ 
This way  we obtain the $C,M,Y$ (cyan - magenta - yellow) colour scheme. These colour systems are mathematically equivalent, but to White should correspond  Black  as the inverse Colour. But Black  does not physically exist in the electro-magnetic spectrum as Colour (in the $RGB$ scheme, Black means an absence of energy). Therefore  in praxis (e.g. in the print industry), Black is artificially added to the $CMY$ system to have $CMYK$ system (the character $K$ is added as the abbreviation for Black). 
\end{remark}

 Subsuming the previous two sections, the following lemma holds:

\begin{lemma} The triple $(\triangle, \oplus, \mathfrak{O})$ is an Abel additive group with respect to Cancellation law.
\end{lemma}

%%%%%%%%%%%%%%%%%%%%%%%%%%%%%%%%%%%%%%%%%%%%%%%%%%%
\section{Multiplication in $\triangle$ \newline \textit{New transformations of Colours} }
%%%%%%%%%%%%%%%%%%%%%%%%%%%%%%%%%%%%%%%%%%%%%%%%%%%%
For simple mixing of two Colours, it is enough to deal with the additive group of Colours. However, there are 
theoretical transformations of Colours which can be called as  multiplications according  their mathematical properties. The author did not know  about any appearance of these operations in the praxis. However, using a computer  digitalization, this theory enables, we can  artificially produce and explore these  Colours. 

This section is about multiplication  of Colours in the Colour space $\triangle$ and about division  in the factorized Colour space  $\triangle | \mathfrak{S}$ where  $\mathfrak{S}$ is the  ideal of singular elements for division  of Colours, cf. Section~\ref{greyideal}.  

%%%%%%%%%%%%%%%%%%%%%%%%%%%%%
\subsection{Cyclic compositions of polar operators in $\mathbb{C}$}\label{squares}
%%%%%%%%%%%%%%%%%%%%%%%%%%%%%
It is easy to see that for the number 3 ($R, G, B$), there are  possible six symmetric Latin squares which respectively yield 6 commutative operations ("multiplications")  $\otimes_i:\mathbb{A} \times \mathbb{A} \to \mathbb{A}, i= 1,2,\dots, 6$ of Colours.  

$$ \begin{array}{c}
%%%%%%%%%%%%%
 \begin{array}{c|ccc}
\otimes_1 &R&G&B\\ \hline
R&R&G&B\\
G&G&B&R\\
B&B&R&G
\end{array}, \hskip5mm
%%%%%%%%%%%%%
\begin{array}{c|ccc}
\otimes_2 &R&G&B\\ \hline
R&G&B&R\\
G&B&R&G\\
B&R&G&B
\end{array}, \hskip5mm
%%%%%%%%%%%%%
\begin{array}{c|ccc}
\otimes_3 &R&G&B\\ \hline
R&B&R&G\\
G&R&G&B\\
B&G&B&R
\end{array}, 
%%%%%%%%%%%%%%%%%%%%%%%%%%%%%%%
\\  \\
%%%%%%%%%%%%%%%%%%%%%%%%%%%%%%%
 \begin{array}{c|ccc}
\otimes_4 &R&G&B\\ \hline
R&R&B&G\\
G&B&G&R\\
B&G&R&B
\end{array}, \hskip5mm
%%%%%%%%%%%%%%%%%%
 \begin{array}{c|ccc}
\otimes_5 &R&G&B\\ \hline
R&B&G&R\\
G&G&R&B\\
B&R&B&G
\end{array}, \hskip5mm
%%%%%%%%%%%%%%%%%%%%%
\begin{array}{c|ccc}
\otimes_6 &R&G&B\\ \hline
R&G&R&B\\
G&R&B&G\\
B&B&G&R
\end{array}.
\end{array}
$$

In this paper, we will deal only with commutative operations $\otimes$ which have the property
$\exists Y \in \mathbb{A}; \forall X \in \mathbb{A} \mid Y \otimes X = X$. Such are cases  
$\otimes_1$, $\otimes_2$, $\otimes_3$. For the reason of cyclic change, we will  only deal with  the Latin square 
$\otimes = \otimes_1$. 

So, let us define the compositions  of poles $\otimes:\mathbb{A} \otimes \mathbb{A} \to \mathbb{A}$  with the following Latin square table.
$$ \begin{array}{c|ccc}
\otimes &R&G&B\\ \hline
R&R&G&B\\
G&G&B&R\\
B&B&R&G
\end{array}
$$

The operation of multiplication $\odot$ in  $\triangle$ is defined in accord both with the table of the composition $\otimes $  and the multiplication in the semi field $T$. Namely,
$$ \mathbf{x} \odot \mathbf{y} = \left(\sum_{i=1,2,3} A ^{(i)}[x_i + \xi_i(f)\varepsilon]\right) \odot \left(\sum_{j =1,2,3} A ^{(j)}[y_j + \chi_j(f)\varepsilon]\right)$$
$$:= \sum_{i=1,2,3}\sum_{j=1,2,3}\left[ A ^{(i)}\otimes A ^{(j)} \right] \left\{[x_i + \xi_i(f)\varepsilon] \cdot [y_j + \chi_j(f)\varepsilon] \right\} $$
\begin{equation} = \sum_{i=1,2,3}\sum_{j=1,2,3}\left[A ^{(i)}\otimes A ^{(j)}\right] [x_i y_j  + \{x_i \chi_j(f) + y_j\xi_i(f)\}\varepsilon],\label{star}\end{equation} 
where $A ^{(i)}, A ^{(j)} \in \mathbb{A}$, $i,j =1,2,3$; $\mathbf{x} = R [x_i + \xi_i(f)\varepsilon] + G[x_2 + \xi_2(f)\varepsilon] + B[x_3 +\xi_3(f)\varepsilon] \in \triangle$; $\mathbf{y} =  R[y_1+\chi_1(f)\varepsilon] + G [y_2 + \chi_2(f)\varepsilon] + G [y_3 + \chi_3(f)\varepsilon] \in \triangle$.
Note that multiplication in Equation~(\ref{star}) is parabolic complex.
  
The proof of the following lemma is evident.

\begin{lemma} \ 

\begin{enumerate} \item The result of operation of multiplication is with respect to  Cancellation law, i.e.,
$$\mathbf{x \odot y} = (\mathbf{x} \oplus \lambda_1) \odot (\mathbf{y} \oplus \lambda_2) \oplus \lambda_3,$$ 
for every $\lambda_1, \lambda_2, \lambda_3 \in \mathfrak{O}$.

\item The  element $$ \mathbf{1} := R[1 + 0(f)\varepsilon] + G[0 + 0(f)\varepsilon] + B[0 +0(f)\varepsilon] \in \triangle$$ 
 is an \textit{unit element} for the operation of multiplication  in $\triangle$ 
 (and hence also $\mathbf{1} + \lambda$, $\lambda \in \mathfrak{O}$, with respect to the congruence given with   Cancellation law).
\end{enumerate}\end{lemma}

%%%%%%%%%%%%%%%%%%%%%%%%%%%%%%%%%%%%%%%%%%%%%%%%
\subsection{Conjugation in $\triangle$ \textit{(Polarization of the light)}}
%%%%%%%%%%%%%%%%%%%%%%%%%%%%%%%%%%%%%%%%%%%%%%%%%%

To define an operation of division, we introduce an operation of \textit{conjugation}. 
\begin{remark} \rm Conjugation physically means  a polarization of light; according to the chosen Latin square $\otimes_1$, a projection will be done to the $R$-axis. 
\end{remark} 

\begin{definition} \rm 
 \sl Let $\mathbf{x} = R[r +  \rho(f)\varepsilon]  +  G[g +  \gamma(f)\varepsilon] +  B[b+  \beta(f)\varepsilon] \in \triangle$. 
 
 We say that an element 
 $\overline{\mathbf{x}} \in \triangle$ is a \textit{conjugation} of the element $\mathbf{x}$ if 
$$\overline{\mathbf{x}} := R[r +  \rho(f)\varepsilon]+ G [b+  \beta(f)\varepsilon] + B [g +  \gamma(f)\varepsilon]  \in \triangle.$$
\end{definition}

\begin{theorem}\label{4.3}
Let $\mathbf{x} \in \triangle$ be as in previous definition, let 
$$ \mathbf{y} = [R(u + \sigma(f)\varepsilon)  + G(v + \chi(f)\varepsilon)] + B(t + \xi(f)\varepsilon) \in \triangle.$$

Then
\begin{enumerate}
\item $\overline{\overline{\mathbf{x}}} = \mathbf{x} \in \triangle$,
\item $\overline{\mathbf{x}\oplus \mathbf{y}} = \overline{\mathbf{x}} \oplus \overline{\mathbf{y}} \in \triangle,$
\item $\overline{\mathbf{x}\odot \mathbf{y}} = \overline{\mathbf{x}} \odot \overline{\mathbf{y}} \in \triangle$,
\item 
$ \mathbf{y}\odot \overline{\mathbf{y}} =  R \Theta$,
where
$$\begin{array}{cl}  \Theta =   & \Big[ \frac{ ( u - v)^2 + (v- t)^2 + (t - u)^2}{2} \\ 
 & + \{(u-v)(\sigma(f) - \xi(f))  \\
& + (v-t)(\xi(f)- \chi(f)) + (t-u)(\chi(f) - \xi(f))\}\varepsilon \Big] \in T.
\end{array}$$ 
\end{enumerate}
\end{theorem}

\proof
The proofs of items 1., 2., 3. are exercises in algebra, we let them to the reader. We prove the last statement 4.
We have:
$$\begin{array}{rcl}
\mathbf{y}  \odot  \overline {\mathbf{y}} & = & \phantom{ \odot }
\{R[u + \sigma(f)\varepsilon] + G[v + \chi(f)\varepsilon] + B[t + \xi(f)\varepsilon]\} \\ 
&& \odot
\{R[u + \sigma(f)\varepsilon] + G[t + \xi(f)\varepsilon] + B[v + \chi(f)\varepsilon]\} = \end{array}$$
using the composition table of poles and the distributive law, we continue:
$$= R\big\{[u + \sigma(f)\varepsilon]^2 + [v + \chi(f)\varepsilon]^2 + [t + \xi(f)\varepsilon]^2\big\}$$
$$+
G\big\{[u + \sigma(f)\varepsilon][v + \chi(f)\varepsilon] + [u + \sigma(f)\varepsilon][t + \xi(f)\varepsilon] 
+ [v + \chi(f)\varepsilon] [t + \xi(f)\varepsilon]\big\}$$
$$ + B\big\{ [u + \sigma(f)\varepsilon][v + \chi(f)\varepsilon] +  [u + \sigma(f)\varepsilon] [t + \xi(f)\varepsilon] +  [v + \chi(f)\varepsilon][t + \xi(f)\varepsilon]\big\}=$$
 By the definition of subtraction,
$$=R \Big[
\{[u + \sigma(f)\varepsilon]^2 + [v + \chi(f)\varepsilon]^2 + [t + \xi(f)\varepsilon]^2\}$$ 
$$- \{[u + \sigma(f)\varepsilon][v + \chi(f)\varepsilon]$$ $$ + [u + \sigma(f)\varepsilon] [t + \xi(f)\varepsilon]$$
$$ + [v + \chi(f)\varepsilon] [t + \xi(f)\varepsilon]\}
\Big] =$$
Now, we have to show that this is a $R$-polarized element. Indeed, we continue:
\begin{equation} \begin{array}{c}
 = R \Big\{[u + \sigma(f)\varepsilon]^2 + [v + \chi(f)\varepsilon]^2 + [t + \xi(f)\varepsilon] ^2 
\\ - [u + \sigma(f)\varepsilon] [v + \chi(f)\varepsilon] 
\\ - [v + \chi(f)\varepsilon] [t + \xi(f)\varepsilon] 
\\ - [t + \xi(f)\varepsilon] [u + \sigma(f)\varepsilon]\Big\} =\label{sucet} \end{array}\end{equation}

$$= R \Big\{ \left[\frac{1}{2}[u + \sigma(f)\varepsilon] ^2 - [u + \sigma(f)\varepsilon] [v + \chi(f)\varepsilon]
 + \frac{1}{2}[v + \chi(f)\varepsilon] ^2\right]$$ 
 $$  + \left[\frac{1}{2} [u + \sigma(f)\varepsilon] ^2 - [u + \sigma(f)\varepsilon] [t + \xi(f)\varepsilon]+ \frac{1}{2} [t + \xi(f)\varepsilon] ^2 \right]$$ 
 $$ + \left[\frac{1}{2} [v + \chi(f)\varepsilon] ^2 - [v + \chi(f)\varepsilon] [t + \xi(f)\varepsilon] + \frac{1}{2} [t + \xi(f)\varepsilon] ^2 \right] \Big\} $$
\begin{equation} \label{==}\begin{array} {cl} = \frac {R}{2} \{ & \Big[ [u + \sigma(f)\varepsilon] -[v + \chi(f)\varepsilon] \big]^2 \\  & + \big[  [u + \sigma(f)\varepsilon] - [t + \xi(f)\varepsilon] \big]^2 \\  & +\big[ [v + \chi(f)\varepsilon] - [t + \xi(f)\varepsilon] \big]^2 \Big\} \end{array} \end{equation}
\begin{equation}
\begin{array}{rl} =  \frac{R}{2} & 
\Big\{ [\{u - v\}+ \{\sigma(f) - \chi(f)\}\varepsilon]^2 + 
[\{v  - t\} + \{ \chi(f) - \xi(f)\}\varepsilon ]^2 \\ 
& + 
[\{t - u \} +\{\xi(f) - \sigma(f))\}\varepsilon ]^2 \Big\}
\end {array}\end{equation}
\begin{equation}\label{H}\begin{array}{cl}  =  R & \Big[\frac{( u - v)^2 + (v- t)^2 + (t - u)^2}{2} \\ 
 & +   \{(u-v)(\sigma(f) - \xi(f))  \\
& + (v-t)(\xi(f)- \chi(f)) + (t-u)(\chi(f) - \xi(f))\}\varepsilon \Big],
\end{array}\end{equation}
Since $( u - v)^2 + (v- t)^2 + (t - u)^2 \geq 0$,  then $\mathbf{y}\odot \overline{\mathbf{y}} =  R \Theta$, where $\Theta \in T$. \qed

%%%%%%%%%%%%%%%%%%%%%%%%%%%%%%%%%%%%%%%
\subsection{The ideal $\mathfrak{S}$\label{greyideal}}
%%%%%%%%%%%%%%%%%%%%%%%%%%%%%%%%%%%%%%
We proved that $\mathbf{y} \odot \overline{\mathbf{y}}$ is an $R$-polarized element in $\triangle$. But what about   elements $\mathbf{y} \in \triangle$ such that $\mathbf{y} \odot \overline{\mathbf{y}} \in \mathfrak{O}$ and $\mathbf{y} \not\in \mathfrak{O}$?  

\begin{definition}\rm
An element $\mathbf{y} \in \triangle$ such that $\mathbf{y} \odot \overline{\mathbf{y}} = \mathfrak{O}$ and 
$\mathbf{y} \neq \mathfrak{O}$ is called a \textsl{singular element}. The set of all singular elements (including  elements of $\mathfrak{O}$) is denoted by $\mathfrak{H}$.  \end{definition}

\begin{lemma}\label{5.5} $$\mathfrak{H} = \mathfrak{S}.$$ \end{lemma}

\proof 
Let $\mathbf{y} = R [u + \sigma(f)\varepsilon] + G [v + \chi(f)\varepsilon] + B [t + \xi(f)\varepsilon]\in \triangle$. 
From Equation~(\ref{H}) it follows that $\mathbf{y}\odot \overline{\mathbf{y}} \in \mathfrak{O}$ if and only if $u = v = t$  
for arbitrary real functions $\sigma(f), \chi(f), \xi(f) \in \mathbb{R}^{[0, + \infty)}$, 
 and every  $u = v = t \geq 0$. \qed

\begin{corollary}Since the result $t=u=v$ of the previous Lemma proof is symmetrical and the analogical result can be obtained using with the cyclic change for the Latin squares  $(\otimes_1), \otimes_2, \otimes_3$,
the ideal $\mathfrak{S}$ of all singular elements for divisions derived from matrices $\otimes_1, \otimes_2, \otimes_3$ is the same.      
\end{corollary}

\begin{theorem}
$\mathfrak{S}$ is a two-sided ideal in the ring $(\triangle, \oplus, \odot, \mathfrak{O})$ with respect to  operation of multiplication $\odot$ with respect to  Cancellation law.
\end{theorem}
\proof
We have to proof:
$$ \mathfrak{O} \varsubsetneqq \mathfrak{S} \varsubsetneqq \triangle,$$
$$\mathbf{x} \in \mathfrak{S}\  \&  \ \mathbf{y} \in \mathfrak{S} \implies  \mathbf{x} + \mathbf{y}  \in \mathfrak{S},$$
and 
$$\mathbf{x} \in \mathfrak{S} \ \&  \ \mathbf{y} \in \triangle \implies  \mathbf{x} \odot \mathbf{y}  \in \mathfrak{S}.$$

The first two assertions are evident. Prove the third assertion.

Let $\mathbf{x} =  R[a +  \rho(f)\varepsilon] + G[a +  \gamma(f)\varepsilon] + B[a+  \beta(f)\varepsilon] \in \mathfrak{S}$ and let $\mathbf{y} = R [u + \sigma(f)\varepsilon] + G [v + \chi(f)\varepsilon] + B [t + \xi(f)\varepsilon] \in \triangle$.

We have: 
$$\mathbf{x} \odot  \mathbf{y}  = 
\{R[a +  \rho(f)\varepsilon] + G[a +  \gamma(f)\varepsilon] + B[a+  \beta(f)\varepsilon] \}$$ 
$$
\odot
\{R [u + \sigma(f)\varepsilon] + G [v + \chi(f)\varepsilon] + B [t + \xi(f)\varepsilon]\} = $$
since $\mathbf{x} \in \mathfrak{S}$,
$$= \{R\rho(f)\varepsilon + G\gamma(f)\varepsilon + B\beta(f)\varepsilon \} 
\odot
 \{R [u + \sigma(f)\varepsilon] + G [v + \chi(f)\varepsilon] + B [t + \xi(f)\varepsilon]\} = $$
since $\varepsilon^2=0$,
$$= \varepsilon \{ R\rho(f) + G\gamma(f) + B\beta(f) \} 
\odot
 \{R u  + G v  + B t \} = $$
using the composition table of poles, 
$$= \varepsilon \{ R [u\rho(f)  + v \beta(f) + t \gamma(f)] + G[v\rho(f) + u\gamma(f) +t\beta(f)] +B[t\rho(f) + v\gamma(f) + u\beta(f)]\}.$$
 We obtained an element in $\mathfrak{S}$.  \qed 

%%%%%%%%%%%%%%%%%%%%%%%%%%%%%%
\subsection{Division in the factor space $\triangle | \mathfrak{S}$ }
%%%%%%%%%%%%%%%%%%%%%%%%%%%%%%

Let $\mathbf{x,y} \in \triangle$ and $\mathbf{y} \notin \mathfrak{S}$.

\begin{remark}\rm
The unit element $\mathbf {1} = R[1+0(f)\varepsilon] +  G[0+0(f)\varepsilon] + B[0+0(f)\varepsilon]$ is  a  mathematical abstraction. In the real world, there is no tristimulus of this kind. However, using this  theoretical object $\mathbf {1}$, we are able to divide real Colours.  
\end{remark}

Division is defined  as follows:
$$\mathbf{x} \oslash \mathbf{y}
:=(\mathbf{x} \odot \mathbf{1}) \oslash \mathbf{y} = \mathbf{x} \odot (\mathbf{1} \oslash \mathbf{y}).$$

Find the element $ \mathbf{1} \oslash \mathbf{y} \in \triangle$ if it exists.

\begin{theorem}\label{4.7} Let $\mathbf{y} = R[u + \sigma(f)\varepsilon] + G[v + \chi(f)\varepsilon] +  B[t + \xi(f)\varepsilon] \not\in \mathfrak{S}$.

 Then
$$\mathbf{1} \oslash \mathbf{y}
= \{ R [u + \sigma(f)\varepsilon] + G[t + \xi(f)\varepsilon] +  B[v + \chi(f)\varepsilon]\} \cdot (2/\Theta), $$
where $\Theta$ is defined in (\ref{H}).
\end{theorem}

\proof
By Theorem~\ref{4.3}, the item 4.,  and Lemma~\ref{5.5},
$$\mathbf{1} \oslash \mathbf{y} = \frac{\overline{\mathbf{y}}}
{\mathbf{y} \odot \overline{\mathbf{y}}}. \ \ \Box $$

From Theorem~\ref{4.7} it follows that we can divide every colour by another colour but elements in the ideal         $\mathfrak{S}$.

%%%%%%%%%%%%%%%%%%%%%%%%%%%%%%%%%%%%%%%%%%%%%%%%%%%%%%%%%%%%%%%%%%%%%%%%%%%%%%%%%%%%%%%%%%%%%%%%%%%%%
\subsection{Compatibility of the additive and multiplicative structures  of $\triangle$}
%%%%%%%%%%%%%%%%%%%%%%%%%%%%%%%%%%%%%%%%%%%%%%%%%%%%%%%%%%%%%%%%%%%%%%%%%%%%%%%%%%%%%%%%%%%%%%%%%
To be complete in  verifying of the axioms of the field, we bring the following more-less evident lemma.
\begin{lemma}\

\begin{enumerate}
\item $\mathbf{1} \not\in \mathfrak{S}$;
\item  $\mathbf{1} \oslash \mathbf{1} = \mathbf{1}$; 
\item Let $\lambda \in \mathfrak{O}$. Then
$\mathbf{1} \neq \lambda$; 
\item
 Let $\mathbf{x,y,z} \in \triangle$. Then
 $$\mathbf{x} \odot (\mathbf{y} \oplus \mathbf{z}) = (\mathbf{x} \odot \mathbf{y}) \oplus (\mathbf{y} \odot \mathbf{z}).$$
 \end{enumerate}
\end{lemma}
\proof (1)(2) The first and second items are trivial.

(3)Let $[a +  \alpha(f)\varepsilon] \in T$. Then  
$$ \mathbf{1} = R[(1+a) + \alpha(f)\varepsilon] + G[a +  \alpha(f)\varepsilon] + B[a +  \alpha(f)\varepsilon] $$ 
$$ \neq R[a +  \alpha(f)\varepsilon] + G[a +  \alpha(f)\varepsilon] + G [a +  \alpha(f)\varepsilon] \in \mathfrak{O}.
$$

(4) The fourth item (distributive law) is ensured by construction of the additive $\oplus$ and multiplicative $\odot$ operations in the set  $\triangle$. The  verification of this equation is elementary and we let it to the reader as an exercise.
\qed

%%%%%%%%%%%%%%%%%%%%%%%%%%%%%%%%%%%%%%
\section{Mathematical subsuming}
%%%%%%%%%%%%%%%%%%%%%%%%%%%%%%%%%%%%%%

We collect mathematical results of the paper  into the following theorem. 

\begin{theorem}  Let $\mathbb{A}\subset \mathbb{C}$ be a set of three poles. Let $T$ be a semi-field of triangular coefficients, cf. Subsection~\ref{coef}. For a fixed Latin square $\otimes\in \{\otimes_1, \otimes_2, \otimes_3\}$, cf. Section~\ref{squares}, the system $(\triangle$, $\oplus$, $\odot$, $\mathfrak{O}$, $\mathbf{1})$ (called the Colour space) is an commutative Abel ring  (with respect to Cancellation law congruence) and particular division, c.f. Subsection~\ref{pooles}. 

There are two Abel groups:  an additive group $(\triangle, \oplus, \mathfrak{O})$ (with respect to Cancellation law). 

The second group $(\triangle | \mathfrak{S}, \odot, \mathbf{1} )$ (with respect to Cancellation law)  is a multiplicative  group with the unit $\mathbf{1}$, where the set $\mathfrak{S}$, see Lemma~\ref{5.5}, is an  ideal, cf.~Section~\ref{achroma}. 

The additive and multiplicative groups are linked together by the distributive law of addition with respect to the multiplication/division  commutatively from both sides.
 
The structure $(\triangle | \mathfrak{S}, \oplus, \odot, \mathfrak{S}, \mathbf{1})$ is an field (with respect to Cancellation law). 
 \end{theorem}

%%%%%%%%%%%%%%%%%%%%%%%%%%%%%%%%%%%%%%%%%%%%%%%%%%%%
\section{Conclusions}
%%%%%%%%%%%%%%%%%%%%%%%%%%%%%%%%%%%%%%%%%%%%%%%%%%%%%%%%%%%%%%%%%%%%%%%%%%%%%%%%%%%%%%%%%%%%%%%%%%%%%%%%%%%%%%%%%%
We created and described  a mathematical theory of Colour space  which factorized with the ideal $\mathfrak{S}$ is a tripolar $RGB$ field with respect to Cancellation law. This theory  has  practical and theoretical application to everything where the phenomenon Colour is sofisticated on the   $RGB$ language. For practical applications, the model needs to include a more thorough  theory of the achromatic part of Colour and also it supposed to take into account some corrections implied from the  technical equipment limitations  and human sensory distortions.

%%%%%%%%%%%%%%%%%%%%%%%%%%%%%%%%%%%%%%%%%%%%%%%%%%%

\end{document}